\documentclass[a4paper,11pt]{article}

\usepackage{amssymb}
\usepackage{amsmath}
\usepackage{amsfonts}
\usepackage{amssymb}
\usepackage{amsthm}
\usepackage{amsfonts}

\newcommand {\Div}{\mathrm{div}}
\newcommand {\R}{\mathbb{R}}
\newcommand {\Rn}{\mathbb{R}^d}
\newcommand{\eps}{\varepsilon}

\newcommand{\ve}{\mathbf{v}}

\newcommand{\no}{\mathbf{n}}

\newcommand{\bvarphi}{\boldsymbol{\varphi}}

\newcommand{\bset}[1]{\bigl\{#1\bigr\}}  %geschweifte Klammern, nächste Größe zu {}
\newcommand{\bigarg}[1]{\left(#1\right)}  %runde Klammern
\newcommand{\bigsqb}[1]{\left[#1\right]}  %eckige Klammern

\newtheorem{theorem}{Theorem}
\newtheorem{remark}{Remark}
\newtheorem{lemma}{Lemma}

\begin{document}

\begin{titlepage}
\title{(Non-)Convergence of Solutions of the Convective Allen-Cahn Equation
}

\date{\today \\[2ex]
\emph{Dedicated to Prof.\ Hideo Kozono on the occasion of his 60th birthday.}}

\author{Helmut Abels\footnote{Faculty of Mathematics, University of Regensburg, 93040 Regensburg, Germany, e-mail: helmut.abels@ur.de
  }}

\end{titlepage}

\maketitle

\begin{abstract}
  We consider the sharp interface limit of a convective Allen-Cahn equation, which can be part of a Navier-Stokes/Allen-Cahn system, for different scalings of the mobility $m_\eps=m_0\eps^\theta$ as $\eps\to 0$. In the case $\theta>2$ we show a (non-)convergence result in the sense that the concentrations converge to the solution of a transport equation, but they do not behave like a rescaled optimal profile in normal direction to the interface as in the case $\theta=0$. Moreover, we show that an associated mean curvature functional does not converge the corresponding functional for the sharp interface. Finally, we discuss the convergence in the case $\theta=0,1$ by the method of formally matched asymptotics.\\[1ex]
  \textbf{Key words:} {Two-phase flow, diffuse interface model,  Allen-Cahn equation, sharp interface limit}\\[1ex]
  \textbf{AMS-Classification:} 76T99, 35Q30,
35Q35,
35R35,
76D05,
76D45 
\end{abstract}

\section{Introduction}
\label{intro}

In this contribution we consider the so-called sharp interface limit, i.e., the limit $\eps\to 0$, of the convective Allen-Cahn equation
\begin{align}
  \partial_t c^\eps + \ve \cdot \nabla c^\eps & = m_\eps\left( \Delta c^\eps -\eps^{-2} f(c^\eps)\right) && \mbox{in }\Omega \times (0,T), \label{simpdiffuse1}\\
  c^\eps|_{\partial\Omega} & = -1 && \mbox{on } \partial \Omega \times (0,T), \label{simpdiffuse3}\\
  \left.c^\eps \right|_{t=0}&= c^\eps_0 &&\mbox{in } \Omega . \label{simpdiffuse4}
\end{align}
Here  $\ve\colon \Omega \times [0,T) \to \mathbb{R}^d$ is a given smooth divergence free velocity field with $\no \cdot \ve|_{\partial\Omega}=0$ and  $c^\eps \colon \Omega \times [0,T) \to \mathbb{R}$ is an order parameter, which will be close to the ``pure states'' $\pm 1$ for small $\eps>0$. Here $f=F'$, where $F\colon\R\to \R$ is a suitable double well potential with globale minima $\pm 1$,  e.g.\ $F(c)=(1-c^2)^2$. $c^\eps$ can describe the concentration difference of two different phases in the case of phase transitions, where the total mass of each phase is not necessarily conserved. 
Moreover,
$\Omega\subseteq\Rn$ is assumed to be a bounded domain with smooth boundary, $m_\eps$ is a (constant) mobility coefficient and $\eps>0$ is a parameter that is proportional to the ``thickness'' of the diffuse interface $\{x\in\Omega : |c^\eps(x,t)|<1-\delta\}$ for $\delta \in (0,1)$.  

The convective Allen-Cahn equation \eqref{simpdiffuse1} is part of the following diffuse interface model for the two-phase flow of two incompressible Newtonian, partly miscible fluids with phase transition
\begin{align}
\partial_t \ve^\eps + \ve^\eps \cdot \nabla \ve^\eps - \Div(\nu(c^\eps) D\ve^\eps) + \nabla p^\eps & = -\eps \Div(\nabla c^\eps \otimes \nabla c^\eps), \label{NSAC1}\\
\Div\, \ve^\eps & = 0,  \\\label{NSAC3}
\partial_t c^\eps + \ve^\eps \cdot \nabla c^\eps & = m_\eps\left( \Delta c^\eps -\eps^{-2} f(c^\eps)\right)
\end{align}
in $\Omega\times (0,T)$,
where $\ve^\eps\colon \Omega\times [0,T)\to \Rn$ is the velocity of the mixture, $D \ve^\eps = \frac{1}{2} (\nabla \ve^\eps + (\nabla \ve^\eps)^T)$, $p^\eps\colon \Omega \times [0,T)\to \R$ is the pressure, and $\nu(c^\eps)>0$ is the viscosity of the mixture. This model can be considered as a model for a two-phase flow with phase transition or an approximation of a classical sharp interface model for a two-phase flow of incompressible fluids with surface tension. Here the densities of the two separate fluids are assumed to be the same. A derivation of this model in a more general form with variable densities can be found in Jiang et al.~\cite{VariableDensityNSAC}. We refer to Gal and Grasselli~\cite{GalGrasselliDCDS} for the existence of weak solutions and results on the longtime behavior of solutions for this model and to Giorgini et al.~\cite{GiorginiGrasselliWuNSAC} for analytic results for a volume preserving variant with different densities.
Mathematically, this system arises if one replaces the Cahn-Hilliard equation in the well-known ``model H'', cf.\ e.g.\ \cite{ModelH,GurtinTwoPhase}, by an Allen-Cahn equation.

With the aid of formally matched asymptotic expansions one can formally show that solutions of this system converge to solutions of the following free boundary value problem %, A.\ et al.\ \cite{AbGaGr} showed that the solutions to the ``model H'' converge of the solutions to the following sharp interface models
\begin{align}\label{eq:Limit1a}
\partial_t \ve + \ve \cdot \nabla \ve - \Div(\nu^\pm D\ve) + \nabla p & = 0 && \mbox{in } \Omega^\pm(t), t\in (0,T)\,,  \\ 
\Div\, \ve & = 0 && \mbox{in } \Omega^\pm(t), t\in (0,T) \,,\\
\left[\ve\right]_{\Gamma_t} & = 0 && \mbox{on } \Gamma_t, t\in (0,T)\,,\\
- \left[ \no_{\Gamma_t} \cdot (\nu^\pm D\ve - p \mathrm{Id})\right]_{\Gamma_t} & = \sigma H_{\Gamma_t}\no_{\Gamma_t} && \mbox{on } \Gamma_t, t\in (0,T) \,,\\
V_{\Gamma_t} - \no_{\Gamma_t} \cdot \ve & = m_0 H_{\Gamma_t} && \mbox{on } \Gamma_t, t\in(0,T)\,,\label{eq:Limit1z}
\end{align}
when $m_\eps=m_0>0$ and 
\begin{align}\label{eq:Limit2a}
\partial_t \ve + \ve \cdot \nabla \ve - \Div(\nu^\pm D\ve) + \nabla p & = 0 && \mbox{in } \Omega^\pm(t), t\in (0,T)  \,, \\
\Div \ve & = 0 && \mbox{in } \Omega^\pm(t), t\in (0,T) \,,\\
\left[\ve\right]_{\Gamma_t} & = 0 && \mbox{on } \Gamma_t, t\in(0,T) \,,\\\label{YoungLaplace}
- \left[ \no_{\Gamma_t} \cdot (\nu^\pm D\ve - p \mathrm{Id})\right]_{\Gamma_t} & = \sigma H_{\Gamma_t} \no_{\Gamma_t} && \mbox{on } \Gamma_t, t\in (0,T) \,, \\\label{eq:Limit2z}
V_{\Gamma_t} - \no_{\Gamma_t} \cdot \ve & = 0 && \mbox{on }\Gamma_t, t\in (0,T) \,,
\end{align}
when $m_\eps=m_0\eps$, $m_0>0$. We will discuss this formal result in the appendix in more detail, cf.~Remark~\ref{rem:2} below. Here $\nu^\pm>0$ are viscosity constants, $\Omega^\pm(t) \subset \Omega$ are open and disjoint such that $\partial \Omega^-(t) = \Gamma_t = \partial \Omega^+(t) \cap \Omega$, $\no_{\Gamma_t}$ denotes the outer normal of $\partial\Omega^-(t)$ and  the normal velocity and the mean curvature
of $\Gamma_t$ are denoted by $V_{\Gamma_t}$ and $H_{\Gamma_t}$, respectively, taken with respect to $\no_{\Gamma_t}$. Furthermore, $\left[\,.\, \right]_{\Gamma_t} $ denotes the jump of a
quantity across the interface in the direction of $\no_{\Gamma_t}$, i.e., $ \left[ f\right]_{\Gamma_t}  (x) = \lim_{h \to 0}(f(x + h \no_{\Gamma_t}) - f(x - h \no_{\Gamma_t}))$ for $x \in \Gamma_t$.

In the case $\nu^+ = \nu^-$ and that the Navier-Stokes equation is replaced by a (quasi-stationary) Stokes system Liu and the author proved rigorously  in \cite{AbelsLiuARMA} that the convergence holds true in the first case $m_\eps=m_0>0$ for sufficiently small times and for well-prepared initial data.
More precisely, it was shown that in a neighborhood of $\Gamma_t$
\begin{equation}\label{eq:Shape}
  c_\eps(x,t)= \theta_0\left(\frac{d_{\Gamma_t}(x)-\eps h_\eps(x,t)}\eps\right) + \mathcal{O}(\eps) 
\end{equation}
(even with $\mathcal{O}(\eps^2)$), where $d_{\Gamma_t}$ is the signed distance function to $\Gamma_t$ and $h_\eps$ are correction terms, which are uniformly bounded in $\eps\in (0,1)$, and $\theta_0\colon \R\to \R$ is the so-called optimal profile that is determined by
\begin{align}
- \theta_0'' + f(\theta_0) & = 0 \mbox{ in } \mathbb{R}, & \theta_0(0) & = 0, & \lim_{z \to \pm\infty}\theta_0(z) & = \pm 1 \,.\label{theta0} 
\end{align} 
This form is important in order to obtain in the limit $\eps\to 0$ the Young-Laplace law \eqref{YoungLaplace}, cf.\ e.g.\ \cite[Section~4]{AbelsGarckeGruen2}.

It is the goal of the present contribution to show that in the case $m_\eps=m_0 \eps^\theta$ with $\theta>2$ the solutions of the convective Allen-Cahn equation \eqref{simpdiffuse1}-\eqref{simpdiffuse3} do not have the form \eqref{eq:Shape} in general. Moreover, we will show that the functional
\[
\left\langle H^\eps , \bvarphi \right\rangle: = \eps \int_\Omega{\nabla c^\eps \otimes \nabla c^\eps : \nabla \bvarphi \, dx} % \quad \bvarphi \in C^\infty_{0,\sigma}(\Omega)
\]
does not converge to the mean curvature functional
 \begin{equation}
   \label{eq:MCF}
 2 \sigma \int_{\Gamma_t}{ \no_{\Gamma_t} \otimes \no_{\Gamma_t} : \nabla \bvarphi \, d \mathcal{H}^{d-1}} = - 2 \sigma \int_{\Gamma_t}{H_{\Gamma_t} \no_{\Gamma_t} \cdot \bvarphi \, d \mathcal{H}^{d-1}} 
\end{equation}
 for all $\bvarphi \in C^\infty_{0,\sigma}(\Omega)=\left\lbrace f \in C^\infty_0(\Omega)^d : \Div f = 0 \right\rbrace$, where
 \[
   \sigma = \frac{1}{2} \int_\mathbb{R}{\left( \theta'_0(z) \right)^2 dz}.
 \]
We note that $H^\eps$ is the weak formulation of the right-hand side of \eqref{NSAC1}, which should converge to a weak formulation of the right-hand side of \eqref{YoungLaplace}. Therefore there is no hope that solutions of the full system \eqref{NSAC1}-\eqref{NSAC3} converge to solutions of the corresponding limit system with \eqref{YoungLaplace} as $\eps\to 0$ in the case that $m_\eps= m_0 \eps^\theta$, $\theta>2$. We note that this effect was first observed for the corresponding Navier-Stokes/Cahn-Hilliard system by Schaubeck and the author in \cite{NonconvergenceConvCH} in the case $\theta>3$. These results are also contained in the PhD-thesis of Schaubeck~\cite{PromotionStefan}. It is not difficult to show that $\left\langle H^\eps , \bvarphi \right\rangle$ converges to \eqref{eq:MCF} if \eqref{eq:Shape} holds true in a sufficiently strong sense. Moreover, in the case $\theta>3$ non-convergence of the Navier-Stokes/Cahn-Hilliard system in the case of radial symmetry and an inflow boundary condition was shown by Lengeler and the author in \cite[Section~4]{AbelsLengeler}. We note that the latter counter example can be adapted to the present case of a Navier-Stokes/Allen-Cahn equation in the case $\theta>2$. %%% Checken!!! 

The structure of this contribution is as follows:
In Section~\ref{secnot} we summarize some preliminaries and notation. Afterwards we prove the nonconvergence result in Section~\ref{secnoncon}. Finally, in Section~\ref{sec:Asymptotics} we discuss briefly the sharp interface limit of the convective Allen-Cahn equation in the case $m_\eps= m_0 \eps^\theta$ with $\theta =0,1$.

\section{Preliminaries and Notation} \label{secnot}
We denote $a \otimes b = \left( a_i b_j \right)^d_{i,j=1}$ for $a,b \in \mathbb{R}^d$ and $A:B = \sum^d_{i,j=1}{A_{ij} B_{ij}}$ for $A,B \in \mathbb{R}^{d \times d}$. We assume that $\Omega \subset \mathbb{R}^d $ is a bounded domain with smooth boundary $\partial \Omega$. Furthermore, we define $\Omega_T = \Omega \times (0,T)$ and $\partial_T \Omega = \partial\Omega \times (0,T)$ for $T>0$. Moreover, $\no_{\partial\Omega}$ denotes the exterior unit normal on $\partial \Omega$. For a hypersurface $\Gamma_t \subset \Omega$, $t\in [0,T]$, without boundary such that $\Gamma_t = \partial \Omega^-(t)$ for a  domain $\Omega^-(t) \subset\subset \Omega$, the interior domain is denoted by $\Omega^-(t)$ and the exterior domain by $\Omega^+(t) := \Omega \backslash (\Omega^-(t) \cup \Gamma_t)$, i.e., $\Gamma_t$ separates $\Omega$ into an interior and an exterior domain. $\no_{\Gamma_t}$ is the exterior unit normal on $\partial \Omega^-(t)=\Gamma_t$. The mean curvature of $\Gamma_t$ with respect to $\no_{\Gamma_t}$ is denoted by $H_{\Gamma_t}$. In the following $d_{\Gamma_t}$ is the signed distance function to $\Gamma_t$ chosen such that $d_{\Gamma_t} < 0$ in $\Omega^-(t)$ and $d_{\Gamma_t}>0$ in $\Omega^+(t) $. By this convention we obtain $\nabla d_{\Gamma_t} = \no_{\Gamma_t}$ on $\Gamma_t$. Moreover, we define
$$Q^\pm := \left\lbrace (x,t) \in \Omega_T : d(x,t) \gtrless 0 \right\rbrace.$$
The ``double-well'' potential $F: \mathbb{R} \to \mathbb{R}$ is a smooth function taking its global minimum $0$ at $\pm1$. For its derivative $f(c)=F'(c)$ we assume
\begin{align}
f(\pm 1) & = 0 , & f'(\pm1) & > 0, & \int^u_{-1}{f(s)\, ds} = \int^u_1{f(s)\, ds} & > 0  \label{bedingungphi1}
\end{align}
for all $u\in (-1,1)$.  In equation (\ref{simpdiffuse1}) the given velocity field satisfies $\ve \in C^0_b([0,T]; C^4_b(\overline{\Omega}))^d$ with $\Div \,\ve = 0$ and $\ve \cdot \no_{\partial \Omega} = 0$ on $\partial \Omega$ and the mobility constant $m_\eps$ has the form $m_\eps = m_0\eps^\theta$ for some $\theta\geq 0$ and $m_0>0$. In equation (\ref{simpdiffuse4}) we choose the special initial value
\begin{align}
\left.c^\eps\right|_{t=0} & = \zeta\!\left(\tfrac{d_{\Gamma_0}}{\delta}\right) \theta_0\!\left(\tfrac{d_{\Gamma_0}}{\eps}\right) + \left( 1- \zeta\!\left(\tfrac{d_{\Gamma_0}}{\delta}\right) \right) \left(2 \chi_{\left\{d_{\Gamma_0} \geq 0\right\}} -1 \right) & \mbox{in }\Omega\,, \label{initialcepsilon}
\end{align}
where we determine the constant $\delta >0 $ later. Here $\zeta \in C^\infty_0(\mathbb{R})$ is a cut-off function such that
\begin{align}
\zeta(z) & = 1 \mbox{ if } \left| z \right| < \frac{1}{2} , & \zeta(z) & = 0 \mbox{ if } \left| z \right| > 1 , & z \zeta'(z) \leq 0 \mbox{ in } \mathbb{R} \,, \label{cutoffzeta}
\end{align}
and $\theta_0$ is the unique solution to (\ref{theta0}).
This choice of the initial value is natural in view of \eqref{eq:Shape}.

\section{Nonconvergence Result} \label{secnoncon}

Our main result is:

\begin{theorem}\label{thmnegres}
Let $\Omega \subset \mathbb{R}^d $ be a bounded domain with smooth boundary $\partial \Omega$, $\Gamma_0$ a smooth hypersurface such that $\Gamma_0=\partial\Omega_0^-$ for a domain $\Omega_0^-\subset\subset \Omega$ and let $c^\eps$ be the solution to the convective Allen-Cahn equation (\ref{simpdiffuse1})-(\ref{simpdiffuse3}) with initial condition (\ref{initialcepsilon}). Then for every $T>0$ and for all $\bvarphi \in C^\infty([0,T];\mathcal{D}(\Omega)^d)$ we have
\begin{eqnarray*}
\int^T_0{\left\langle H^\eps , \bvarphi \right\rangle dt} \to_{\eps\to 0} 2 \sigma \int^T_0{\int_{\Gamma_t}{ \left| \nabla (d_{\Gamma_0}(X^{-1}_t)) \right| \no_{\Gamma_t} \otimes \no_{\Gamma_t} : \nabla \bvarphi \, d \mathcal{H}^{d-1} } \, dt} \,,
\end{eqnarray*}
where the evolving hypersurface $\Gamma_t$, $t\in [0,T]$, is the solution of the evolution equation
\begin{eqnarray*}
V_{\Gamma_t}(x) = \no_{\Gamma_t}(x,t) \cdot \ve(x,t) \mbox{ for }x\in \Gamma_t, t \in (0,T], \quad \Gamma(0)=\Gamma_{0},
\end{eqnarray*}
where $V_{\Gamma_t}$ is the normal velocity of $\Gamma_t$. Moreover, it holds
\begin{eqnarray*}
\left\| c^\eps -(2 \chi_{Q^+} -1 ) \right\|^2_{L^2(\Omega_T)} = \mathcal{O}(\eps) \quad \text{as }\eps \to 0.
\end{eqnarray*}
\end{theorem}

\begin{remark}
  In general $\left| \nabla (d_{\Gamma_0}(X^{-1}_t)) \right| = \left| D X^{-T}_t \nabla d_{\Gamma_0} \circ X^{-1}_t \right| \neq 1$, we refer to \cite[Remark~1]{NonconvergenceConvCH} for a proof. This shows that the weak formulation of $H^\eps$ does not converge to the weak formulation of the right-hand side of the Young-Laplace law \eqref{YoungLaplace} in general.
\end{remark}

To prove the theorem we follow the same strategy as in \cite{NonconvergenceConvCH}: First we construct a family of approximate solutions $\left\lbrace c^\eps_A\right\rbrace_{0<\eps \leq 1}$. Afterwards we estimate the difference $\nabla (c^\eps -c^\eps_A)$, which will enable us to prove the assertion of the theorem.
We start with the observation that $\Gamma_t:= X_t(\Gamma_0)$ is the solution to the evolution equation.

\begin{lemma} \label{lemmahypersur}
Let $\Gamma_{0} \subset \Omega$ be a given smooth hypersurface such that $\Gamma_0=\partial\Omega_0^-$ for a domain $\Omega_0^-\subset\subset\Omega$. Then the evolving hypersurface $\Gamma_t:= X_t\left(\Gamma_{0}\right) \subset \Omega$, $t\in [0,T]$, is the solution to the problem 
\[
V_{\Gamma_t}= \no_{\Gamma_t} \cdot \ve \quad \mbox{on } \Gamma_t, t \in (0,T), \quad \Gamma(0)=\Gamma_{0}.
\]
\end{lemma}
We refer to \cite[Lemma~3]{NonconvergenceConvCH} for the proof.

For the following let $P_{\Gamma_t}(x)$ be the orthogonal projection of $x$ onto $\Gamma_t$. Then there exists a constant $\delta >0 $ such that $\Gamma_t(\delta):= \left\{ x \in \Omega : \left|d_{\Gamma_t}(x))\right| < \delta \right\} \subset \Omega$ and $\tau_t \colon \Gamma_t(\delta) \to (-\delta,\delta) \times \Gamma_t$ defined by $\tau_t(x) = (d_{\Gamma_t}(x),P_{\Gamma_t}(x))$ is a smooth diffeomorphism, cf.\ e.g.\ \cite[Kapitel 4.6]{HildebrandtAnalysis2}. %%% Pruess Simonett einfügen.
% Furthermore, we define
% \begin{align*}
%   \Gamma &:=  \left\lbrace (x,t) \in \Omega_T : d_{\Gamma_t}(x) = 0 \right\rbrace\qquad \text{and} \\
%  \Gamma(\delta) &: = \left\lbrace (x,t) \in \Omega_T : \left| d_{\Gamma_t}(x) \right| < \delta  \right\rbrace.  
% \end{align*}

We will need the following result:
\begin{lemma} \label{lemmaevolution}
For $e \colon \bigcup_{t \in [0,T]} X_t(\Gamma_0(\delta)) \times \left\{ t \right\} \to \mathbb{R}$ defined by $e(x,t):= d_{\Gamma_0}(X_t^{-1}(x))$ the following properties hold:
\begin{enumerate}
	\item $ \frac{d}{dt} e(x,t) = - \ve(x,t) \cdot \nabla e(x,t) $ for all $(x,t) \in \bigcup_{t \in [0,T]} X_t(\Gamma_0(\delta)) \times \left\{ t \right\} $\,.\label{elemma1}
	\item $ e(x,t)$ is a level set function for $\Gamma_t $, i.e., $e(x,t)=0 $ if and only if $ x \in \Gamma_t$.\label{elemma2}
\end{enumerate}
\end{lemma}
We refer to  \cite[Lemma~4]{NonconvergenceConvCH} for the proof.

As mentioned in Section \ref{secnot}, let $\theta_0$ be the solution to (\ref{theta0}) and let $\zeta$ be a cut-off function as in (\ref{cutoffzeta}). Then we define
\[
c_A^\eps(x,t):= \left\{ \begin{array}{l@{\;}l}
\pm1 & \mbox{in } \overline{Q^\pm} \cap \bigcup\limits_{t \in [0,T]} \overline{X_t(\Omega \backslash \Gamma_0(\delta))} \times \{ t \} \,, \\
\zeta\!\left(\tfrac{e}{\delta}\right) \theta_0\!\left(\tfrac{e}{\eps}\right) \pm (1-\zeta\!\left(\tfrac{e}{\delta}\right)) & \mbox{in } Q^\pm \cap \bigcup\limits_{t \in [0,T]} X_t(\Gamma_0(\delta) \backslash \Gamma_0\!\left(\tfrac{\delta}{2}\right)) \times \{ t \}, \\
\theta_0\!\left(\tfrac{e}{\eps}\right) & \mbox{in } \bigcup\limits_{t \in [0,T]} X_t(\Gamma_0\!\left(\tfrac{\delta}{2}\right)) \times \{t\} \,.
\end{array} \right.
\]
Then we have $c^\eps_A(.,0) = c^\eps(.,0)$ since $e(.,0) = d_{\Gamma_0}$ and
\begin{equation*}
  \partial_t c^\eps_A + \ve \cdot \nabla c^\eps_A = 0\qquad \text{in }\Omega_T  
\end{equation*}
since $\partial_t e + \ve \cdot \nabla e = 0$. Moreover, by the construction
\begin{equation*}
  c^\eps_A|_{\partial\Omega}=0 \qquad \text{on }\partial\Omega.
\end{equation*}
Furthermore, we define the approximate mean curvature functional by 
\[
\left\langle H^\eps_A ,\varphi \right\rangle = \eps \int_\Omega{\nabla c^\eps_A \otimes \nabla c^\eps_A : \nabla \bvarphi \, dx}\,.
\]
for all $\bvarphi \in \mathcal{D}(\Omega)^d$. Then we have:

\begin{lemma} \label{lemmarestterme}
Let $c^\eps_A$ be defined as above. Then there exists some constant $C>0$ independent of $\eps$  and $\eps_0 \in (0,1] $ such that the estimates
\begin{eqnarray}
\left\| \Delta c^\eps_A(.,t) \right\|_{L^2(\Omega)} & \leq & C \eps^{- \frac{3}{2}} \,,\label{DeltaepsilonA}\\\label{NablaepsilonA}
\left\| \nabla c^\eps_A(.,t) \right\|_{L^2(\Omega)} & \leq & C \eps^{- \frac{1}{2}} \,, \\
\left\| f(c^\eps_A(.,t)) \right\|_{L^2(\Omega)} & \leq & C \eps^{\frac{1}{2}} \,, \label{fepsilonA} \\
\left\| c^\eps_A(.,t) - (2 \chi_{Q^+}(.,t) -1) \right\|_{L^2(\Omega)} & \leq & C \eps^{\frac{1}{2}}  \label{cepsilonA+-1}
\end{eqnarray}
hold for all $t \in [0,T]$ and $\eps \in (0,\eps_0)$.
\end{lemma}
We refer to  \cite[Lemma~5]{NonconvergenceConvCH} for the proof. % In the same way as \eqref{fepsilonA} one proves
% \begin{equation*}
%   \left\| f(c^\eps_A(.,t)) \right\|_{L^1(\Omega)}  \leq  C \eps.
% \end{equation*}

Now we are able to prove the central lemma for the proof of Theorem~\ref{thmnegres}.

\begin{lemma} \label{lemmanablaR}
Let $c^A_\eps$ be defined as above and let $c^\eps$ be the unique solution to \eqref{simpdiffuse1}-\eqref{simpdiffuse3} with initial condition \eqref{initialcepsilon}. Then, for $\theta \geq 2$, there exists some constant $C >0 $ independent of $\eps$ and $\eps_0>0$ such that 
\begin{align}
 \eps \left\| \nabla (c^\eps - c^\eps_A) \right\|^2_{L^2(\Omega_T)} &\leq C \eps^{\theta -2} \quad \text{and}\\
\left\| c^\eps - c^\eps_A \right\|_{L^\infty(0,T;L^2(\Omega))} &\leq C \eps^{\theta - \frac32}
\end{align}
for all $\eps \in (0,\eps_0]$.
\end{lemma}
\proof First of all, we note that $c^\eps(x,t),c^\eps_A(x,t)\in [-1,1]$ for all $x\in\Omega$, $t\in (0,T)$. For $c^\eps_A$ this follows from the construction and for $c^\eps$ by the maximum principle.

We denote by $u = c^\eps - c^\eps_A$ the difference between exact and approximate solution, which solves
\begin{equation*}
  \partial_t c^\eps_A + \ve \cdot \nabla c^\eps_A = 0\qquad  \text{in }\Omega_T.  
\end{equation*}
We multiply the difference of the differential equations for $ c^\eps$ and $ c^\eps_{A}$ by $u$ and integrate the resulting equation over $\Omega$. Then we get for all $t \in (0,T)$
\begin{align*}
  0 =&   \int_\Omega u\left[\partial_t u + \ve \cdot \nabla u - m_0\eps^{\theta} \Delta u - m_0\eps^{\theta} \Delta c^\eps_A + m_0\eps^{\theta - 2}  f(c^\eps) u\right] dx \\
  =&   \int_\Omega\left(\partial_t \tfrac{|u|^2}2 - \ve\cdot \nabla \tfrac{|u|^2}2+ m_0 \eps^{\theta} |\nabla u|^2\right) \, dx\\
     & + \int_\Omega\left( m_0\eps^{\theta} \nabla u \cdot \nabla c^\eps_A - m_0\eps^{\theta-2}  f(c^\eps) u\right) \, dx \\
    = &   \frac{1}{2} \frac{d}{dt} \int_\Omega | u|^2 dx + m_0 \eps^\theta\int_\Omega|\nabla u|^2 dx + \int_\Omega\left(m_0 \eps^{\theta} \nabla u \cdot \nabla c^\eps_A + m_0\eps^{\theta -2 } u f(c^\eps)\right) dx \,,
\end{align*} 
where we have used $u = 0$ on $\partial \Omega$ as well as $\Div \, \ve = 0$ in $\Omega$. 
%Using $f(s)  \geq 0$ for $s \geq C_0 \geq 1$ and $f(s) \leq 0$ for $s \leq -C_0 \leq -1$ as well as $|c_A(x,t)|\leq 1$, we have
%\begin{eqnarray*}
%\int_{\left\{x \in \Omega \,:\, \left|c^\eps(x,t)\right|\geq C_0\right\}}{f(c^\eps) u \, dx}\geq 0\,.
%\end{eqnarray*}
By Hölder's and Young's inequalities we obtain
\begin{align}
  & \frac{1}{2} \frac{d}{dt} \int_\Omega{ \left|u \right|^2 dx} + \frac{m_0}2\eps^{\theta} \int_\Omega{\left| \nabla u \right|^2 dx } %+ m_0\eps^{\theta-2}\int_{\left\{x \in \Omega \,:\, \left|c^\eps(x,t)\right|\geq C_0\right\}}{f(c^\eps) u \, dx }
    \nonumber\\
& \leq  C\eps^{\theta}\left\| \nabla c^\eps_A\right\|_{L^2(\Omega)}^2 + \eps^{\theta-2} \left| \int_{\Omega}{f(c^\eps) u \, dx} \right| \label{absch1lemmaR}
\end{align}
for all $\eps \in (0,\eps_0)$, where
\begin{align}
  \left|\int_{\Omega}{f(c^\eps) u \, dx}\right| & \leq  \left|\int_{\Omega}{f(c_A^\eps) u \, dx}\right|  +C\|u\|_{L^2(\Omega)}^2 \nonumber \\
      & \leq  \left\| f(c^\eps_A) \right\|_{L^2(\Omega)}\|u\|_{L^2(\Omega)} + C \left\| u\right\|^2_{L^2(\Omega)}   \nonumber\\
  &\leq  C \eps^{\frac12}\| u\|_{L^2(\Omega)} + C \left\| u\right\|^2_{L^2(\Omega)} \label{fcepsilon}
\end{align}
since $f'$ is Lipschitz continuous on $[-1,1]$.
 Hence (\ref{absch1lemmaR}) together with (\ref{fcepsilon}) and (\ref{NablaepsilonA}) yield 
 \begin{align*}
   & \frac{1}{2} \frac{d}{dt} \int_\Omega{ \left|u \right|^2 dx} + m_0\eps^{\theta} \int_\Omega{\left| \nabla u \right|^2 dx } \nonumber \\
   & \leq  C \left( \left\| u \right\|^2_{L^2(\Omega)} + \eps^{2\theta - 3} +\eps^{\theta-2}\|u\|_{L^2(\Omega)}^2 \right)  \leq  C_1 \left( \left\| u \right\|^2_{L^2(\Omega)} + \eps^{2\theta - 3} \right) 
\end{align*}
since $\theta\geq 2$ for some $C_1>0$ independent of $\eps$ and $t \in [0,T]$.
Hence the Gronwall inequality implies
\begin{eqnarray*}
\sup_{0 \leq t \leq T} \left\|u \right\|^2_{L^2(\Omega)} + \eps^\theta \|\nabla u\|_{L^2((0,T)\times\Omega)}^2 \leq C \eps^{2\theta-3} 
\end{eqnarray*}
for some $C= C(T) >0$ independent of $\eps$. Therefore the lemma is proved.
\hfill $\Box$\\

Now we can show that $H^\eps-H^\eps_A$ converges to $0$ as $\eps$ goes to zero.

\begin{lemma} \label{lemmaKruemmung1}
Let $H^\eps$ and $H^\eps_A$ be defined as above and let $\theta > 2$. Then it holds 
\begin{eqnarray*}
\left| \int_0^T{\left\langle H^\eps - H^\eps_A , \bvarphi \right\rangle dt} \right| \to_{\eps\to 0} 0\,,
\end{eqnarray*}
for all $\bvarphi \in C^\infty([0,T];\mathcal{D}(\Omega)^d)$.
\end{lemma}

\proof The proof is almost the same as in \cite[Lemma~6]{NonconvergenceConvCH}. But we include it for the convenience of the reader since the argument is central for our main result. Let $\bvarphi \in C^\infty([0,T];\mathcal{D}(\Omega)^d)$ and set $u = c^\eps -c^\eps_A$. Then 
\begin{align*}
& \eps \left| \int_{\Omega_T}{\left( \nabla c^\eps \otimes \nabla c^\eps - \nabla c^\eps_A \otimes \nabla c^\eps_A \right) : \nabla \bvarphi \, dx} \right|  \\
& \leq  \eps \left| \int_{\Omega_T}{\left( \nabla c^\eps \otimes \nabla u \right) : \nabla \bvarphi \, dx} \right| + \eps \left| \int_{\Omega_T}{\left( \nabla u \otimes \nabla c^\eps_A \right) : \nabla \bvarphi \, dx} \right| \\
& \leq  \eps \left\|\nabla \bvarphi \right\|_{L^\infty(\Omega_T)} \left\| \nabla u\right\|_{L^2(\Omega_T)} \left( \left\| \nabla c^\eps\right\|_{L^2(\Omega_T)} + \left\| \nabla c^\eps_A\right\|_{L^2(\Omega_T)} \right) .
\end{align*}
Because of Lemma \ref{lemmarestterme} and Lemma~\ref{lemmanablaR}, we have
\[
\left\| \nabla c^\eps \right\|_{L^2(\Omega_T)} \leq \left\| \nabla c^\eps_A \right\|_{L^2(\Omega_T)} + \left\| \nabla u \right\|_{L^2(\Omega_T)}\leq C\left(\eps^{-\frac12}+ \eps^{\frac{\theta-3}2}\right). 
\]
Using Lemma \ref{lemmanablaR} we conclude
\begin{eqnarray*}
\left| \int_0^T{\left\langle H^\eps - H^\eps_A , \bvarphi \right\rangle dt} \right| \leq C \eps^{\frac{ \theta - 2}{2}} \left( 1 + \eps^{\frac{\theta - 2}2} \right) 
\end{eqnarray*}
for some constant $C=C(\varphi)>0$ and for all $\eps$ small enough. Since $\theta >2 $, the assertion follows. \hfill $\Box$\\

\begin{lemma} \label{lemmaKruemmung2}
Let $H^\eps_A$ and $c^\eps_A$ be defined as above. Then it holds for all $\bvarphi \in \mathcal{D}(\Omega)^d$ and $t \in [0,T]$
\begin{align*}
\left\langle H^\eps_A, \bvarphi \right\rangle & \to_{\eps\to 0} 2 \sigma \int_{\Gamma_t}{ \left| \nabla (d_{\Gamma_0}(X^{-1}_t)) \right| \no_{\Gamma_t} \otimes \no_{\Gamma_t} : \nabla \bvarphi \, d \mathcal{H}^{d-1} }.
\end{align*} 
\end{lemma}
We refer to  \cite[Lemma~8]{NonconvergenceConvCH} for the proof.

\medskip

\noindent \textbf{Proof of Theorem \ref{thmnegres}:} The first assertion of the theorem immediately follows by Lemma \ref{lemmaKruemmung1} and \ref{lemmaKruemmung2}. 
The second assertion is a consequence of Lemma \ref{lemmarestterme} and Lemma \ref{lemmanablaR} since $\theta > 3$.
\makebox[1cm]{} \hfill $\Box$

\section{Formal Asymptotics}\label{sec:Asymptotics}

In this section we will use the method of formally matched asymptotic expansions to identify the sharp interface limit of the convective Allen-Cahn equation  \eqref{simpdiffuse1}-\eqref{simpdiffuse3} in the cases $m_\eps= m_0 \eps^\theta$ for $\theta=0,1$ and some $m_0>0$. We follow similar arguments as in \cite[Section~4]{AbelsGarckeGruen2}. 
In particular we assume that there are smoothly evolving hypersurfaces $\Gamma_t$, $t\in (0,T)$, such that $\Gamma_t= \partial\Omega^-(t)$, and we have the following expansions:\\[1ex]
 \emph{Outer expansion:} ``Away from $\Gamma_t$''  we assume that  $c_\eps$ has an expansion of the form:
    \begin{equation*}
      c_\eps(x,t)= \sum_{k=0}^\infty  \eps^k {c_k^\pm(x,t)} \qquad \text{for every }x\in\Omega^\pm(t).
    \end{equation*}
 \emph{Inner expansion:} In a neighborhood $\Gamma_t(\delta)$, $\delta>0$, of $\Gamma_t$ $c_\eps$ has an expansion of the form:
    \begin{equation*}
      c_\eps(x,t)= \sum_{k=0}^\infty  \eps^k {c_k(\tfrac{d_{\Gamma_t}}\eps, P_{\Gamma_t}(x),t)} \quad \text{for all }x\in \Gamma_t(\delta).
    \end{equation*}
  \emph{Matching condition:} 
    \begin{alignat*}{2}
      \lim_{z\to\pm  \infty} c_k(z,x,t) &= c_k^\pm (x,t) &\qquad& \text{for all }x\in\Gamma_t, k=0,1,\\
      \lim_{z\to\pm  \infty} \partial_z c_0(z,x,t) &= 0 &\qquad& \text{for all }x\in\Gamma_t. 
    \end{alignat*}
  Moreover, all functions in the expansions above are assumed to be sufficiently smooth.

In the following we will use the expansions above and the matching conditions, insert them into the convective Allen-Cahn equation \eqref{simpdiffuse1} and equate all  terms of same order in order to determine the leading parts in the inner and outer expansions formally.

	\subsection{Outer Expansion} 
	First we use a power series expansion of $c_\eps$ due to the outer expansion. Then
	\begin{equation*}
	f'(c_\varepsilon(x,t)) = f' (c_0^\pm(x,t)) c_1^\pm (x,t)+ \varepsilon f'' (c_0^\pm (x,t)) c_1^\pm (x,t) + \mathcal{O} (\varepsilon^2)
	\end{equation*}
        and
	we obtain from \eqref{simpdiffuse1}
	\begin{equation*}
	\frac{1}{\varepsilon^{2-k}} f' (c_0^\pm (x,t))+ \frac{1}{\varepsilon^{1-k}} f''(c_0^\pm (x,t))c_1^\pm(x,t)
	+ \mathcal{O}(1) = 0
	\end{equation*}
        for all $x\in \Omega^\pm(t)$. This yields
	\begin{itemize}
        \item [i)] At order $\frac{1}{\eps^{2-k}}$ we obtain $f'(c_0^\pm (x,t)) =0$. Thus $c_0^\pm (x,t) \in \bset{\pm 1,0}$.  Here we exclude the case $c_0^\pm (x,t) =0$ since $0$ is unstable and define $\Omega^\pm (t)$ such that
          \begin{equation*}
            c_0^\pm (x,t) =\pm 1 \text{ for all } x\in \Omega^\pm (t).
          \end{equation*}
	\item[ii)] If $k=0$, we obtain at order $\frac{1}{\eps}$ that $f''(c_0 (x,t)) c_1^\pm (x,t) =0$. Since $f''(\pm 1)>0$, we conclude
          \begin{equation*}
            c_1^\pm (x,t) =0 \text{ for all } x\in \Omega^\pm (t).
          \end{equation*}
          If $k=1$, the corresponding term is of order $\mathcal{O}(1)$ and we do not use this information. Moreover, we will not determine $c^\pm_1$ and $c_1$ in this case.
	\end{itemize}

	\subsection{Inner Expansion}
	In $\Gamma_t(\delta)$ we use the inner expansion in \eqref{simpdiffuse1} in order to determine the leading coefficients $c_0(\rho,s,t)$ and, in the case $k=0$, $c_1(\rho,s,t)$, where $s:=s(x):= P_{\Gamma_t} (x)$. To this end we use
	\begin{alignat*}{1}
          \ve\cdot \nabla c_j (\rho, s,t) &=  \frac{1}{\varepsilon} \ve\cdot \nabla d_{\Gamma_t} (\rho,s,t)  +\mathcal{O}(1),\\
	\Delta c_j (\rho,s,t) &=  \frac{1}{\varepsilon^2} (\partial_\rho^2 c_j) \bigarg{\rho,s,t}
+  \frac{1}{\varepsilon} (\partial_\rho c_j) \bigarg{\rho , s,t}\Delta d_{\Gamma_t} (x)  + \mathcal{O}(1),\\
	\partial_t c_j (\rho ,s,t)&
         =  \frac{1}{\varepsilon} (\partial_\rho c_j) \bigarg{\rho,s,t} \partial_t d_{\Gamma_t} (x) +\mathcal{O}(1)
       \end{alignat*}
       on  $\Gamma_t$,   where $\rho = \tfrac{d_{\Gamma_t} (x,t)}{\varepsilon}$ and
       \begin{equation*}
       \nabla d_{\Gamma_t}= \no_{\Gamma_t},\quad \Delta d_{\Gamma_t} = -H_{\Gamma_t},\quad  \partial_t d_{\Gamma_t} =-V_{\Gamma_t} \qquad \text{on }\Gamma_t.  
       \end{equation*} 
	Hence inserting the inner expansion in \eqref{simpdiffuse1} and equating terms of the same order yields for all $x \in {\Gamma_t}$: 
	\begin{alignat*}{1}
          & m_0\bigsqb{-\partial_\rho^2 c_0 (\rho,s,t) + f'(c_0(\rho,s,t))} \cdot \frac{1}{\varepsilon^2}\\
           & + m_0\left[ -\partial_\rho^2 c_1 (\rho,s,t) + f''(c_0(\rho,s,t)) c_1(\rho,s,t)\right]\cdot \frac1\eps\\
        & + \left[- \partial_\rho c_0 (\rho,s,t) (V_{\Gamma_t} -\no_{\Gamma_t}\cdot \ve - m_0 H_{\Gamma_t} )\right] \cdot \frac{1}{\varepsilon}=  O (1)
      \end{alignat*}
      in the case $k=0$ and
	\begin{alignat*}{1}
          & \bigsqb{m_0\left(-\partial_\rho^2 c_0 (\rho,s,t) + f'(c_0(\rho,s,t))\right)- (\partial_\rho c_0) (\rho,s,t) (V_{\Gamma_t}-\no_{\Gamma_t}\cdot \ve)} \cdot \frac{1}{\varepsilon}	 =  O (1)
      \end{alignat*}
      in the case $k=1$. For the following we distinguish the cases $k=0,1$.

      \noindent
      \emph{Case $k=0$:} The $\mathcal{O}(\frac{1}{\varepsilon^2})$-terms yield
	\begin{equation*}
	 - \partial_\rho^2 c_0 (\rho,s,t) + f'(c_0(\rho,s,t)) =0\qquad \text{for all }\rho\in\R, s\in \Gamma_t, t\in [0,T].
	\end{equation*}
	Because of the matching condition, we obtain
        $$
        \lim_{\rho \to \pm \infty} c_0(\rho,s,t) =c_0^\pm (s,t) = \pm 1 \qquad \text{for all }s\in \Gamma_t, t\in [0,T].
        $$
	In order to obtain that $\Gamma_t$ approximates the zero-level set of $c_\eps(x,t)= c_0(\tfrac{d_{\Gamma_t}}\eps,s(x),t)+\mathcal{O}(\eps)$ sufficiently well, we obtain $c_0(0,s,t) = 0$. %%% Anders schreibe???
        Hence
	\begin{alignat*}{1}
	c_0(\rho,x,t) & = \theta_0 (\rho) \quad \text{for all } x \in \Gamma_t, \rho\in\R.
	\end{alignat*}
        Furthermore, the $\mathcal{O}(\frac{1}{\varepsilon})$-terms yield
        \begin{equation*}
          m_0\left(-\partial_\rho^2 c_1 (\rho,x,t) + f''(\theta_0 (\rho)) c_1 (\rho,x,t)\right)  = \theta'_0(\rho) (V_{\Gamma_t} - \no_{\Gamma_t}\cdot \ve - m_0H_{\Gamma_t})=:g(\rho)
        \end{equation*}
	Since $\theta'_0$ is in the kernel of the differential operator $-\partial_\rho^2 +f''(\theta_0)$, this ODE has a bounded solution if and only if
	\begin{equation}\label{eq:CompCond1}
	\int_{\R} g(\rho) \theta'_0 (\rho) d \rho =0, 
	\end{equation}
        which  is equivalent to
        \begin{equation*}
          V_{\Gamma_t} - \no_{\Gamma_t}\cdot \ve = H_{{\Gamma_t} }\qquad \text{ on } {\Gamma_t}. 
        \end{equation*}
	Now the matching condition yields $c_1 (\rho,x,t) \to_{\rho \to \pm \infty} c_1^\pm \equiv 0$. Hence
	$c_1 \equiv 0$ since the solution is unique. Altogether we obtain for the inner expansion
	\begin{equation*}
	c_\varepsilon (x,t) = \theta_0 \bigarg{\frac{d_{\Gamma_t} (x)}{\varepsilon}} + O (\varepsilon^2)
	\end{equation*}
	close to $\Gamma_t$.

        \noindent
      \emph{Case $k=1$:} The $\mathcal{O}(\frac{1}{\varepsilon})$-terms yield
	\begin{align}\nonumber
          &m_0\left(- \partial_\rho^2 c_0 (\rho,s,t) + f'(c_0(\rho,s,t))\right)\\\label{eq:k1}
          &\quad -\partial_\rho c_0(\rho,s,t) (V_{\Gamma_t} (s)-\no_{\Gamma_t}(s)\cdot \ve(s,t)) =0
	\end{align}
        for all $s\in\Gamma_t$. Testing with $\partial_\rho c_0(\rho,x,t)$ yields
        \begin{align*}
          0 = \int_{\R}|\partial_\rho c_0(\rho,s,t)|^2\, d\rho \left(V_{\Gamma_t} (s)-\no_{\Gamma_t}\cdot \ve(s,t) \right)
        \end{align*}
        since
        \begin{equation*}
          \int_{\R}\partial_\rho \left(\frac{|\partial_\rho c_0(\rho,s,t)|^2}2 + f(c_0(\rho,s,t))\right) d\rho =0
        \end{equation*}
        because of the matching condition for $\partial_\rho c_0$. Because of $c_0(\rho,s,t)\to_{\rho\to \pm\infty} \pm 1$,  $\partial_\rho c_0$ does not vanish and we obtain
        \begin{equation*}
          V_{\Gamma_t} = \no_{\Gamma_t}\cdot \ve \qquad \text{ on } {\Gamma_t}. 
        \end{equation*}
        Moreover, we obtain from \eqref{eq:k1}
        \begin{equation*}
          - \partial_\rho^2 c_0 (\rho,s,t) + f'(c_0(\rho,s,t))= 0 \qquad \text{for all }s\in \Gamma_t,\rho\in\R.
        \end{equation*}
        Hence we can conclude as in the case $k=0$ that $c_0(\rho,s,t)= \theta_0(\rho)$ for all $\rho\in\R$ and $s\in\Gamma_t$, $t\in [0,T]$.
        \begin{remark}\label{rem:2}
          The formal calculations show that $c_\eps$ should have an expansion of the form \eqref{eq:Shape} in the case $\theta=0,1$. This is important to obtain \eqref{YoungLaplace} in the limit.
          Actually, using $c_0(\rho,s,t)=\theta_0(\rho)$ one can easily modify the results in \cite[Section~4]{AbelsGarckeGruen2} to show formally convergence of the Navier-Stokes/Allen-Cahn system \eqref{NSAC1}-\eqref{NSAC3} to \eqref{eq:Limit1a}-\eqref{eq:Limit1z} in the case $\theta=0$ and \eqref{eq:Limit2a}-\eqref{eq:Limit2z} in the case $\theta=1$. A rigorous justification of this convergence under suitable assumptions remains open.
        \end{remark}

%\bibliographystyle{spmpsci}      % mathematics and physical sciences
%\bibliography{../../Bibliography.bib}   % name your BibTeX data base

\def\cprime{$'$} \def\ocirc#1{\ifmmode\setbox0=\hbox{$#1$}\dimen0=\ht0
  \advance\dimen0 by1pt\rlap{\hbox to\wd0{\hss\raise\dimen0
  \hbox{\hskip.2em$\scriptscriptstyle\circ$}\hss}}#1\else {\accent"17 #1}\fi}

 \end{document}